\documentclass[reqno]{amsart}
\usepackage{amsmath,amssymb,amsthm,mathtools}
\usepackage{microtype}
\usepackage{enumitem}
\usepackage{needspace}
\usepackage{cite}
\usepackage[hidelinks]{hyperref}
\hypersetup{pdftitle={Every finite incline is finitely based}}
\renewcommand{\le}{\leqslant}
\renewcommand{\ge}{\geqslant}
\theoremstyle{plain}
\newtheorem{theorem}{Theorem}[section]
\newtheorem{lemma}[theorem]{Lemma}
\newtheorem{proposition}[theorem]{Proposition}
\newtheorem{corollary}[theorem]{Corollary}
\theoremstyle{definition}

\newtheorem{remark}[theorem]{Remark}
\newcommand{\AI}{\Sigma_{\mathrm{ai}}}
\newcommand{\W}{\mathcal W}
\newcommand{\V}{\mathcal V}

\newcommand{\var}{\operatorname{var}}
\newcommand{\con}{\operatorname{con}}
\newcommand{\Mod}{\operatorname{Mod}}
\newcommand{\Cminus}[1]{\mathbf C_{#1}^{-}}
\newcommand{\doi}[1]{\href{https://doi.org/#1}{\nolinkurl{#1}}}
\setlist[enumerate]{label=\textup{(\roman*)},leftmargin=2.1em,itemsep=3pt}
\title{Every finite incline is finitely based}

\author{Haihong Wu}
\address{School of Mathematics, Northwest University, Xi'an, 710127, Shaanxi, P.R. China}
\email{17836198299@163.com}

\date{}
\makeatletter
\@namedef{subjclassname@2020}{\textup{2020} Mathematics Subject Classification}
\makeatother


\begin{document}
\begin{abstract}
We prove that every incline of finite additive height generates a Cross
variety. In particular, every finite incline is finitely based. Moreover,
we obtain a uniform bound on the number of variables needed to axiomatize
subvarieties of decreasing Catalan semiring varieties, and give an
effective construction of finite equational bases.
\end{abstract}
\maketitle

\section{Introduction}

An \emph{additively idempotent semiring} (or \emph{ai-semiring} for short)
is an algebra $(S;+,{\cdot})$ such that $(S;+)$ is a commutative
idempotent semigroup, $(S;{\cdot})$ is a semigroup, and multiplication
distributes over addition on both sides. The relation $a\le b$ defined
by $a+b=b$ is a partial order compatible with both operations. An
ai-semiring is called an \emph{incline} if it satisfies
\begin{equation}\label{eq:contraction}
 xy\le x,\qquad xy\le y.
\end{equation}
Multiplication in an incline need not be commutative. Inclines were
studied systematically by Cao, Kim, and Roush~\cite{CaoKimRoush1984};
see also the survey of Kim and Roush~\cite{KimRoush2004} and the
representation theorem of Kim, Roush, and
Markowsky~\cite{KimRoushMarkowsky1997}. Related classes of
semilattice-ordered semigroups and additively idempotent semirings are
considered in~\cite{KurilPolak2005,KepkaKorbelar2025}.

A class of algebras of the same signature is a \emph{variety} if it is
closed under taking subalgebras, homomorphic images, and arbitrary direct
products. By Birkhoff's theorem, varieties are precisely the classes
defined by identities. For an algebra $S$, let $\var S$ denote the
variety generated by $S$. The algebra $S$ and the variety $\var S$
satisfy the same identities. A variety is \emph{finitely based} if it
can be defined by a finite set of identities; otherwise, it is
\emph{nonfinitely based}. An algebra is finitely based if the variety
it generates is finitely based. The \emph{finite basis problem} for a
class of algebras asks which members of that class are finitely based.

The finite basis problem for ai-semirings has received considerable
attention. Dolinka constructed a finite nonfinitely based ai-semiring
and a finite inherently nonfinitely based
ai-semiring~\cite{Dolinka2007,Dolinka2009INFB}. Further examples arise
from finite inverse semigroups~\cite{GusevVolkov2023} and the Brandt
monoid~\cite{Volkov2021}. Jackson, Ren, and Zhao~\cite{JacksonRenZhao2022}
showed that an ai-semiring can be nonfinitely based even when its
multiplicative reduct is finitely based. Semirings related to $S_7$
and the extension $S_7^0$ have been studied
in~\cite{GaoJacksonRenZhao2025,WuRenZhao2024}.

Ai-semirings of small order provide further examples and positive
results. Shao and Ren~\cite{ShaoRen2015} studied the varieties generated
by ai-semirings of order two. Zhao et
al.~\cite{ZhaoRenCrvenkovicShaoDapic2020}, together with the results of
Jackson et al.~\cite{JacksonRenZhao2022}, settled the finite basis
problem for ai-semirings of order three. Recent work on ai-semirings
of order four includes the series
\cite{RenLiuZengChen2025I,YueRenZengShao2025II,
RenLiuYueChen2026III,RenYueYuanLyuYangYu2026IV} and the nonfinite-basis
results of Yue and Ren~\cite{YueRen2026NFB} and Yue, Ren, and
Gao~\cite{YueRenGao2026Two}.

Finite basability is also closely related to the structure of
subvariety lattices. Ren and Zhao~\cite{RenZhao2016}, Ren, Zhao,
and Wang~\cite{RenZhaoWang2017}, and Ren, Zhao, and
Shao~\cite{RenZhaoShao2020} studied varieties subject to power
identities. Hereditary finite basability was investigated by Ren and
Zeng~\cite{RenZeng2019}. On the other hand, ai-semiring varieties can
have uncountable subvariety lattices~\cite{Dolinka2008}, and flat
extensions of groups give rise to limit
varieties~\cite{RenJacksonZhaoLei2023}. A further investigation of
intervals in subvariety lattices appears in the manuscript of Gao et
al.~\cite{GaoJacksonRenZhaoContinuum}.

A variety is \emph{hereditarily finitely based} if every subvariety is
finitely based. It is a \emph{Cross variety} if it is finitely based,
generated by a finite algebra, and has only finitely many subvarieties;
see~\cite{MacdonaldVaughanLee1978}. Every subvariety of a Cross variety
is again Cross. In particular, every Cross variety is hereditarily
finitely based. Our main result is the following.

\begin{theorem}\label{thm:main}
Every finite incline generates a Cross variety.
In particular, every finite incline is finitely based.
\end{theorem}

The proof uses a family of varieties considered by Gusev and
Volkov~\cite{GusevVolkovCatalan}. For $n\ge3$, let $\W_n$ be the
variety of inclines satisfying
\begin{equation}\label{eq:sigma}
 \sigma_n:\quad
 x_1\cdots x_{n-1}\approx
 \sum_{i=1}^{n-1}x_1\cdots x_{i-1}x_i^{\,n}x_{i+1}\cdots x_{n-1},
\end{equation}
and put
\begin{equation}\label{eq:K}
 K_n=n-2+2^{\,n-1}.
\end{equation}
Gusev and Volkov proved that $\W_n=\var\Cminus n$, where
$\Cminus n$ is the decreasing Catalan
semiring~\cite[Theorem~4.12]{GusevVolkovCatalan}. They also showed that
every incline of finite additive height $h$ belongs to
$\var\Cminus{h+2}$~\cite[Corollary~5.9]{GusevVolkovCatalan}. Thus the
identities $\sigma_n$ and the description of $\W_n$ as a generated
variety are known. We prove the following stronger finite-basis result.

\begin{theorem}\label{thm:uniform}
For every $n\ge3$, the variety $\W_n$ is locally finite and has only
finitely many subvarieties. Every subvariety $\V\subseteq\W_n$ is
finitely based and generated by a finite algebra. Moreover, all
identities of $\V$ follow from the defining identities of $\W_n$ and
the identities of $\V$ involving at most $K_n$ variables.
\end{theorem}

In~\cite[Section~6.3, p.~36]{GusevVolkovCatalan}, Gusev and Volkov ask
how many subvarieties $\var\Cminus n$ has for $n\ge3$ and note that
no finite nonfinitely based incline is known. Theorems~\ref{thm:main}
and~\ref{thm:uniform} answer these questions: every finite incline is
finitely based, and each of these decreasing Catalan semiring varieties
has a finite subvariety lattice.

The main step is the uniform bound in Theorem~\ref{thm:uniform}. We
first express every term as a sum of block words with at most $n-2$
marked positions. We then identify variables according to their
occurrences in the idempotent blocks. A second substitution recovers
the left-hand side of an inequality and decreases its right-hand side.
Local finiteness converts this variable bound into a finite basis and
a finite subvariety lattice.

Throughout the paper, we work in the constant-free signature
$\{+,{\cdot}\}$, so all substitution terms are nonempty. We also
sharpen the finite-height reduction: an incline of additive height
$h\ge2$ belongs to $\W_{h+1}$. Consequently,
Theorem~\ref{thm:main} extends to all inclines of finite additive
height, including those with infinite underlying sets.

\section{Preliminaries and the finite-height reduction}

An \emph{additively idempotent semiring}, or \emph{ai-semiring}, is an
algebra $(S;+,\cdot)$ satisfying
\begin{align*}
 x+y&\approx y+x,& (x+y)+z&\approx x+(y+z),&x+x&\approx x,\\
 (xy)z&\approx x(yz),&x(y+z)&\approx xy+xz,&(x+y)z&\approx xz+yz.
\end{align*}
Let $\AI$ denote this set of six identities. The relation
$a\le b$ if and only if $a+b=b$ is a partial order, addition is its
join operation, and both operations preserve the order in each argument.
For example, if $a\le b$, then
$ac+bc=(a+b)c=bc$, so $ac\le bc$; the other cases are similar.
An inequality between terms is always an abbreviation for an identity:
$u\le v$ means $u+v\approx v$.

The two laws in~\eqref{eq:contraction}, together with $\AI$, define
inclines. Let $\Delta_n$ be the finite set consisting of $\AI$,
these two laws written as identities, and $\sigma_n$. Thus
$\W_n=\Mod(\Delta_n)$. When we establish an identity in every member
of $\W_n$, it follows from $\Delta_n$ by completeness of equational
logic. We use this observation throughout without repeating it.

Fix a countably infinite alphabet of variables. A \emph{word} is a nonempty finite
product of variables, and $\con(w)$ is the set of variables occurring
in a word $w$. The same notation is used for terms. By distributivity,
every term is equivalent modulo $\AI$ to a finite nonempty sum of
words. Whenever a displayed product contains an empty block, that block
is omitted. This convention introduces no identity element or constant
symbol.

The \emph{height} of a partially ordered set is the maximum cardinality
of a chain, when this maximum is finite. In particular, height counts
elements, not strict inequalities. A variety is \emph{locally finite}
if each of its finitely generated members is finite. We write
$F_{\V}(d)$ for the free algebra of a variety $\V$ on $d$ generators,
where $d\ge1$. Standard facts about varieties and free algebras may be
found in~\cite{BurrisSankappanavar1981}.

\begin{lemma}\label{lem:deletion}
In an incline, deleting factors from a product can only increase its
value, provided the remaining product is nonempty. In particular,
$a\ge a^2\ge a^3\ge\cdots$.
\end{lemma}

\begin{proof}
It suffices to delete one factor. At an end of the product this follows
from~\eqref{eq:contraction}. For an interior factor, write the product
as $uav$, with $u$ and $v$ nonempty. Since $ua\le u$, monotonicity
gives $uav\le uv$. Repeated deletion proves the assertion.
\end{proof}

\begin{lemma}\label{lem:bottom}
Every incline of finite additive height has a least element $0$ in its
natural order, and $0a=a0=0$ for all $a$.
\end{lemma}

\begin{proof}
Finite height ensures that the order has a minimal element $e$: starting
at any element, a strictly descending sequence must terminate.
For any $a$, contraction gives $ea\le e$ and $ea\le a$.
Minimality of $e$ implies $ea=e$, and hence $e\le a$. Thus $e$ is a
least element; denote it by $0$. Again by contraction,
$0a\le0$ and $a0\le0$, so both products equal $0$.
The symbol $0$ here names an element of the algebra and does not belong
to the signature.
\end{proof}

\begin{proposition}\label{prop:height}
Let $S$ be an incline whose additive semilattice has finite height
$h\ge2$. Then $S$ satisfies $\sigma_{h+1}$, and consequently
$S\in\W_{h+1}$.
\end{proposition}

\begin{proof}
Set $n=h+1$, take arbitrary $a_1,\ldots,a_h\in S$, and write
$p_j=a_1\cdots a_j$. By contraction,
\[
 p_1\ge p_2\ge\cdots\ge p_h.
\]
For each $i$, let
\[
 t_i=a_1\cdots a_{i-1}a_i^{\,n}a_{i+1}\cdots a_h.
\]
Lemma~\ref{lem:deletion} gives $t_i\le p_h$ for every $i$.

Suppose first that $p_{i-1}=p_i$ for some $2\le i\le h$.
Then $p_{i-1}a_i=p_{i-1}$, and induction yields
$p_{i-1}a_i^r=p_{i-1}$ for every positive integer $r$.
It follows that $t_i=p_h$. Thus $\sum_{i=1}^h t_i=p_h$.

Otherwise $p_1>p_2>\cdots>p_h$ is a chain of $h$ distinct elements.
Let $0$ be the least element supplied by Lemma~\ref{lem:bottom}.
If $p_h\ne0$, appending $0$ gives a chain of $h+1$ distinct
elements, contrary to the definition of $h$. Hence $p_h=0$.
Since $t_i\le p_h$, all the $t_i$ equal $0$, and their sum again
equals $p_h$. This proves $\sigma_n$.
\end{proof}

An incline of height one has just one element, since the join of any
two distinct elements would be strictly above at least one of them.
It satisfies every identity. Thus the height-one case can always be
handled separately, or included by using $\W_{\max\{3,h+1\}}$.

\section{Idempotent blocks and normal forms}

Fix $n\ge3$. All identities in this section are considered in $\W_n$.
Substituting $x$ for every variable in $\sigma_n$ and using additive
idempotency gives
\[
 x^{n-1}\approx x^{2n-2}.
\]
By contraction, $x^{2n-2}\le x^n\le x^{n-1}$, and therefore
\begin{equation}\label{eq:powers}
 x^{n-1}\approx x^n.
\end{equation}
Multiplying this identity by further powers of $x$ shows that all powers
with exponent at least $n-1$ coincide. In particular, $t^n$ is
multiplicatively idempotent for every term $t$.

Fix a total order on the variable alphabet. This orders the variable
symbols and is independent of the natural order in an incline. For a
nonempty finite set $A=\{y_1,\ldots,y_s\}$, listed in that order, put
\begin{equation}\label{eq:block}
 E_A=(y_1\cdots y_s)^n.
\end{equation}

\Needspace{10\baselineskip}
\begin{lemma}\label{lem:greatest}
Under every evaluation in a member of $\W_n$, the value of $E_A$ is
the greatest multiplicatively idempotent common lower bound of the
values of the variables in $A$. Consequently, for every word $w$,
\begin{equation}\label{eq:content}
 w^n\approx E_{\con(w)}.
\end{equation}
In particular, the value of $E_A$ does not depend on the chosen order of
the variables in $A$.
\end{lemma}

\begin{proof}
Let $S\in\W_n$, and evaluate $y_j$ at $a_j\in S$ for
$1\le j\le s$. The element $(a_1\cdots a_s)^n$ is multiplicatively
idempotent by~\eqref{eq:powers} and is below every $a_j$ by
Lemma~\ref{lem:deletion}. Now let $e\in S$ be any element such that
$e^2=e$ and $e\le a_j$ for all $j$. Monotonicity gives
\[
 e=e^{sn}\le(a_1\cdots a_s)^n.
\]
Thus $(a_1\cdots a_s)^n$ is the greatest element among the
multiplicatively idempotent common lower bounds of $a_1,\ldots,a_s$.
The same argument applies to the factors of $w$, including repeated
factors. Repetitions do not change the set of common lower bounds.
Hence $w^n$ and $E_{\con(w)}$ have the same greatest-element
characterization and are equal under every evaluation. This also
proves independence of the chosen order.
\end{proof}

\begin{lemma}\label{lem:absorb}
In any incline, if $e^2=e$ and $e\le a$, then $ea=ae=e$.
\end{lemma}

\begin{proof}
Monotonicity and contraction give $e=e^2\le ea\le e$ and
$e=e^2\le ae\le e$.
\end{proof}

A \emph{block word} is an expression
\begin{equation}\label{eq:blockword}
 E_{A_0}a_1E_{A_1}\cdots a_mE_{A_m},
\end{equation}
where each $a_i$ is an individual variable and each $A_i$ is a finite
set of variables. Empty blocks are omitted. The $m$ individual letter
positions occupied by $a_1,\ldots,a_m$ are called \emph{marked
positions}. Marked letters may repeat, and may occur in the blocks as
well. We allow $m=0$, in which case $A_0$ must be nonempty.
The wholly empty expression is never a block word.

\begin{proposition}\label{prop:normal}
Modulo $\Delta_n$, every word is a finite nonempty sum of block words
with at most $n-2$ marked positions. The same assertion holds for every
semiring term.
\end{proposition}

\begin{proof}
Every word is a block word with all blocks empty. It therefore
suffices to prove the assertion for block words. We proceed by induction
on the number $m$ of marked positions. If $m\le n-2$, the block word
itself is a sum with one summand of the required form.

Suppose that $m\ge n-1$, and put $k=n-1$. Write $E_i=E_{A_i}$ when
$A_i\ne\varnothing$, omitting $E_i$ otherwise. For $1\le i\le k$,
put
\[
 b_i=E_{i-1}a_iE_i.
\]
Each $b_i$ is a legitimate nonempty substitution word, since it contains
$a_i$. Each intermediate block occurs twice in $b_1\cdots b_k$;
as those blocks are idempotent, we obtain
\begin{equation}\label{eq:duplicated}
 b_1\cdots b_k\approx E_0a_1E_1\cdots a_kE_k.
\end{equation}
Substituting $b_i$ for $x_i$ in $\sigma_n$ gives
\[
 b_1\cdots b_k\approx
 \sum_{i=1}^k b_1\cdots b_{i-1}b_i^n b_{i+1}\cdots b_k.
\]
If $m>k$, multiply both sides on the right by
$a_{k+1}E_{k+1}\cdots a_mE_m$. By~\eqref{eq:duplicated}, the
left-hand side is then the original block word. If $m=k$, no suffix
is needed.

In the $i$th summand on the right, only $b_i$ has been replaced by
$b_i^n$.
By Lemma~\ref{lem:greatest},
\[
 b_i^n\approx F_i,
 \qquad F_i=E_{A_{i-1}\cup\{a_i\}\cup A_i}.
\]
If $E_{i-1}$ is present, then $F_i\le E_{i-1}$: the value of $F_i$
is an idempotent lower bound for every variable in $A_{i-1}$, so
Lemma~\ref{lem:greatest} applies. Similarly, $F_i\le E_i$ whenever
$E_i$ is present. Lemma~\ref{lem:absorb} therefore gives
\begin{equation}\label{eq:neighbours}
 E_{i-1}F_i=F_i=F_iE_i
\end{equation}
whenever the indicated neighboring block exists.

Replace $b_i^n$ by $F_i$, remove the other duplicated blocks by
idempotency, and absorb any adjacent copies of $E_{i-1}$ and $E_i$
using~\eqref{eq:neighbours}. The resulting summand has the form of
\eqref{eq:blockword}, with the segment $E_{i-1}a_iE_i$ replaced by
the single block $F_i$. This replacement results from applying
$\sigma_n$ and reducing the powered factor $b_i^n$. The letter $a_i$ is now
contained in an idempotent block rather than occupying a marked
position. Hence the number of marked positions decreases by one.
The same calculation applies at $i=1$ and $i=k$, with any missing
exterior factors omitted, and when either neighboring block is empty.
Each summand therefore has $m-1$ marked positions, so the induction
hypothesis completes the reduction.

Finally, expand an arbitrary term into a finite nonempty sum of words
by distributivity, and apply the preceding reduction to each word.
\end{proof}

\begin{corollary}\label{cor:local}
The variety $\W_n$ is locally finite. More precisely, for every $d\ge1$,
put
\begin{equation}\label{eq:B}
 B_n(d)=\sum_{m=0}^{n-2}d^m2^{d(m+1)}-1.
\end{equation}
Then
\begin{equation}\label{eq:freebound}
 |F_{\W_n}(d)|\le 2^{B_n(d)}-1.
\end{equation}
\end{corollary}

\begin{proof}
On a fixed alphabet of $d$ variables, a block word with $m$ marked
positions is specified by an ordered list of $m$ marked letters and
$m+1$ subsets $A_0,\ldots,A_m$ of the alphabet. There are $d^m$
choices for the marked letters and $2^d$ choices for each subset.
Thus there are $d^m2^{d(m+1)}$ specifications.
The subtraction of one in~\eqref{eq:B} excludes the wholly empty
expression for $m=0$. Different specifications may represent the same
word, so $B_n(d)$ is an upper bound for the number of block words
needed in Proposition~\ref{prop:normal}.

By Proposition~\ref{prop:normal}, every term is equivalent to a
nonempty sum of these block words. Commutativity and idempotency of
addition allow us to ignore their order and remove repetitions, so
there are at most $2^{B_n(d)}-1$ possible sums. The elements of
$F_{\W_n}(d)$ are the resulting equivalence classes of terms, which
proves~\eqref{eq:freebound}. Every $d$-generated member of $\W_n$ is a
homomorphic image of this finite free algebra and is therefore finite.
\end{proof}

\section{A uniform bound on the number of variables}

We now prove the main step. In the following lemma the upper term $q$
need not have been reduced to normal form.

\begin{lemma}\label{lem:compression}
Let $\V\subseteq\W_n$, and suppose that $\V$ satisfies an inequality
\begin{equation}\label{eq:target}
 v\le q,\qquad
 v=E_{A_0}a_1E_{A_1}\cdots a_mE_{A_m},\quad m\le n-2.
\end{equation}
Then~\eqref{eq:target} follows from $\Delta_n$ and a single identity of
$\V$ involving at most $m+2^{m+1}$ variables.
\end{lemma}

\begin{proof}
Let $B=\{a_1,\ldots,a_m\}$ and
$X=\con(v)\cup\con(q)$. Thus $|B|\le m$.
For $x,y\in X\setminus B$, define
\begin{equation}\label{eq:patterns}
 x\sim y\quad\Longleftrightarrow\quad
 (x\in A_i\Longleftrightarrow y\in A_i)
 \text{ for every }0\le i\le m.
\end{equation}
For each variable in $X\setminus B$, record one bit for each of
the $m+1$ sets $A_0,\ldots,A_m$: the bit is $1$ if the variable
belongs to that set and $0$ otherwise. By~\eqref{eq:patterns}, two
variables are equivalent precisely when these binary vectors agree.
There are at most $2^{m+1}$ equivalence classes, regardless of
their sizes.
The all-zero vector is included whenever it occurs, so variables
appearing only in $q$ are retained.

\Needspace{10\baselineskip}
For each equivalence class $G$, choose a fresh variable $z_G$. Define a substitution
$\varphi$ that fixes the variables in $B$ and sends every variable in
$G$ to $z_G$. Then
\begin{equation}\label{eq:compressed}
 \varphi(v)\le\varphi(q)
\end{equation}
holds in $\V$, being a substitution instance of~\eqref{eq:target}.
Its variables consist of those in $B$ and one fresh variable for
each equivalence class. It therefore involves at most
$|B|+2^{m+1}\le m+2^{m+1}$ variables.

Define a second substitution $\psi$ by fixing $B$ and setting
$\psi(z_G)=E_G$ for each equivalence class $G$. Every equivalence class is nonempty, so all
these substitution terms are permitted in our signature. Let
$\theta=\psi\circ\varphi$ on $X$. If $x\in B$, then $\theta(x)=x$;
if $x\in G$, then $\theta(x)=E_G\le x$ by contraction. Thus
$\theta(x)\le x$ for every $x\in X$, in particular for every
variable occurring in $q$. Monotonicity of the term operations gives
\begin{equation}\label{eq:upper}
 \Delta_n\models\theta(q)\le q.
\end{equation}

We claim that $\theta$ fixes every nonempty block of $v$ modulo
$\Delta_n$. Fix $A_i\ne\varnothing$ and an evaluation in a member of
$\W_n$. In the following comparison, terms denote their values under
this evaluation. Put $e=E_{A_i}$. Since $\theta(x)\le x$ for all
$x\in A_i$, monotonicity gives $\theta(E_{A_i})\le E_{A_i}$.

For the reverse inequality, $e$ is idempotent and below the value of
each variable in $A_i$. If $x\in A_i\cap B$, then $\theta(x)=x$, so
$e\le\theta(x)$. If $x\in A_i\setminus B$ belongs to $G$, the
definition of the membership patterns implies $G\subseteq A_i$.
Thus $e$ is an idempotent lower bound for every variable in $G$,
and Lemma~\ref{lem:greatest} gives $e\le E_G=\theta(x)$.

Write $A_i=\{y_1,\ldots,y_s\}$ in the fixed order. We have just shown
that $e\le\theta(y_j)$ for each $j$. Hence
\[
 e=e^{sn}\le
 \bigl(\theta(y_1)\cdots\theta(y_s)\bigr)^n
 =\theta(E_{A_i}).
\]
It follows that $\theta(E_{A_i})\approx E_{A_i}$ holds in $\W_n$.
All marked letters are fixed as well, and therefore
\begin{equation}\label{eq:lower}
 \Delta_n\models\theta(v)\approx v.
\end{equation}

Apply $\psi$ to~\eqref{eq:compressed}. Combining the result
with~\eqref{eq:upper} and~\eqref{eq:lower} yields
\[
 v\approx\theta(v)\le\theta(q)\le q.
\]
The auxiliary identities hold throughout $\W_n$ and so follow from
$\Delta_n$. Thus~\eqref{eq:target} follows from $\Delta_n$ and the
single identity~\eqref{eq:compressed}, as required.
\end{proof}

\begin{proposition}\label{prop:variables}
For every subvariety $\V\subseteq\W_n$, all identities of $\V$
follow from $\Delta_n$ and the identities of $\V$ involving at most
$K_n=n-2+2^{n-1}$ variables.
\end{proposition}

\begin{proof}
An identity $p\approx q$ is equivalent to the two inequalities
$p\le q$ and $q\le p$. Consider the first of these. By
Proposition~\ref{prop:normal},
$\Delta_n\models p\approx v_1+\cdots+v_t$, where every $v_j$ is a
block word with at most $n-2$ marked positions. Since $v_j\le p$
holds in $\W_n$, the inequality $v_j\le q$ holds in $\V$.
Lemma~\ref{lem:compression} derives it using an identity with at most
\[
 m+2^{m+1}\le n-2+2^{n-1}=K_n
\]
variables. Adding the inequalities for $j=1,\ldots,t$ gives $p\le q$.
The same argument proves $q\le p$.
\end{proof}

\begin{remark}\label{rem:noncomm}
The compression argument uses the order characterization of the
idempotent blocks, not a commutativity law for multiplication.
In particular, it does not require the product of arbitrary
multiplicative idempotents to be idempotent. The only absorption used
is Lemma~\ref{lem:absorb}, with its explicit order hypothesis.
\end{remark}

\section{Finite bases, finite generators, and subvarieties}

\begin{proof}[Proof of Theorem~\ref{thm:uniform}]
Local finiteness is Corollary~\ref{cor:local}, and the variable bound is
Proposition~\ref{prop:variables}. Set $K=K_n$ and choose a finite set
$R$ of representatives for all terms on $x_1,\ldots,x_K$ modulo the
identities of $\W_n$. Such a set exists because $F_{\W_n}(K)$ is
finite. For $\V\subseteq\W_n$, put
\begin{equation}\label{eq:Gamma}
 \Gamma_{\V}=\{r\approx s:r,s\in R,\ \V\models r\approx s\}.
\end{equation}
This is a finite set of identities of $\V$.

Every identity of $\V$ involving at most $K$ variables can first be
renamed onto $x_1,\ldots,x_K$. Replace both of its terms by their
representatives in $R$. These replacements follow from $\Delta_n$,
and the resulting identity belongs to $\Gamma_{\V}$. Renaming back
shows that all such identities follow from
\begin{equation}\label{eq:finitebasis}
 \Delta_n\cup\Gamma_{\V}.
\end{equation}
By Proposition~\ref{prop:variables}, this finite set is a basis for
all identities of $\V$.

The set $\Gamma_{\V}$ is contained in the finite set of all pairs of
terms in $R$, and~\eqref{eq:finitebasis} determines $\V$ uniquely.
Consequently, $\W_n$ has only finitely many subvarieties.

It remains to prove that each $\V$ has a finite generator. Let
$F=F_{\V}(K)$. This algebra is finite by local finiteness. Since
$F\in\V$, we have $\var F\subseteq\V$. The varieties $\V$ and
$\var F$ have the same identities on at most $K$ variables. Indeed,
one inclusion of their sets of such identities follows from
$F\in\V$. Conversely, if an identity on $x_1,\ldots,x_K$ holds in
$F$, evaluate it at the corresponding free generators. Equality there
means, by the construction of the free algebra, that the identity holds
in $\V$.

Both varieties lie in $\W_n$. Applying
Proposition~\ref{prop:variables} to $\var F$ shows that all its
identities follow from $\Delta_n$ and its identities on at most $K$
variables. All those identities hold in $\V$, so every identity of
$\var F$ holds in $\V$. Hence $\V=\var F$.
\end{proof}

\Needspace{7\baselineskip}
\begin{corollary}\label{cor:hereditary}
For every $n\ge3$, every subvariety of $\W_n$ is a Cross variety.
In particular, $\W_n$ is a Cross variety.
\end{corollary}

\begin{proof}
Theorem~\ref{thm:uniform} supplies a finite basis and a finite
generator for every subvariety. Its own subvarieties form a subset of
the finite subvariety lattice of $\W_n$.
\end{proof}

\begin{corollary}\label{cor:height}
Every incline of finite additive height generates a Cross
variety. If its height is $h\ge2$, the identities of every subvariety
of its generated variety follow from $\Delta_{h+1}$ and identities
involving at most $h-1+2^h$ variables.
\end{corollary}

\begin{proof}
For $h\ge2$, Proposition~\ref{prop:height} places the incline in
$\W_{h+1}$, so apply Theorem~\ref{thm:uniform} and
Corollary~\ref{cor:hereditary}. For $h=1$, the algebra has one element
and generates the trivial variety, which is Cross.
\end{proof}

\begin{proof}[Proof of Theorem~\ref{thm:main}]
The additive semilattice of a finite incline has finite height.
Corollary~\ref{cor:height} therefore applies.
\end{proof}

The proof also gives explicit bounds. Let
$N_n=2^{B_n(K_n)}-1$, with $B_n$ as in~\eqref{eq:B}.
One may choose $|R|\le N_n$. Hence
\begin{equation}\label{eq:latticebound}
 |L(\W_n)|\le2^{N_n^2},
\end{equation}
where $L(\W_n)$ is the subvariety lattice. Moreover, every
$\V\subseteq\W_n$ is generated by an algebra with at most $N_n$
elements. No claim of optimality is made for these bounds.

\section{Effective bases and decreasing Catalan semirings}

\subsection{A finite construction from operation tables}

We now give a finite construction of an equational basis from the
operation tables of a finite incline. The construction uses the normal
form in Proposition~\ref{prop:normal} and the variable bound in
Proposition~\ref{prop:variables}.

Let $S$ be a finite incline of additive height $h\ge2$, and set
\[
 n=h+1,\qquad K=K_n=h-1+2^h,\qquad X=\{x_1,\ldots,x_K\}.
\]
Let $\mathcal B_{n,K}$ be the finite set of all nonempty block words
on $X$ with at most $n-2$ marked positions, using the increasing
order of indices to define each block $E_A$. Let
$\mathcal P_{n,K}$ consist of all nonempty sums of distinct members of
$\mathcal B_{n,K}$, with the summands listed in a fixed order.
Set
\begin{equation}\label{eq:algorithm}
 \Sigma(S)=\Delta_n\cup
 \{P\approx Q:P,Q\in\mathcal P_{n,K},\ S\models P\approx Q\}.
\end{equation}

\begin{proposition}\label{prop:effective}
The set $\Sigma(S)$ is a finite identity basis for $S$, and it can be
constructed effectively from the two operation tables of $S$.
\end{proposition}

\begin{proof}
The set is finite by construction. Every identity in $\Delta_n$ holds
in $S$ by Proposition~\ref{prop:height}, and all remaining identities
were selected because they hold in $S$.

For completeness, take an identity of $S$ involving at most $K$
variables and rename it onto $X$. By Proposition~\ref{prop:normal},
each side is equivalent modulo $\Delta_n$ to a term in
$\mathcal P_{n,K}$. The identity between the resulting terms is
included in~\eqref{eq:algorithm}. Thus $\Sigma(S)$ proves every
identity of $S$ on at most $K$ variables.
Proposition~\ref{prop:variables}, applied to $\var S\subseteq\W_n$,
now gives all identities of $S$.

The addition table determines the order and hence its height $h$.
The two finite sets of terms can be enumerated. For any pair $P,Q$,
the test $S\models P\approx Q$ requires only evaluation on all
$|S|^K$ assignments to $X$, using the operation tables.
All steps therefore terminate.
\end{proof}

For the one-element incline, a basis is $\AI\cup\{x\approx y\}$.
If the height is not computed, one may instead take $n=|S|+2$.
Indeed, the descending sequence of $n-1$ prefix products must then
contain an adjacent equality, and the proof of
Proposition~\ref{prop:height} gives $\sigma_n$ directly.
This proves that a finite basis can be computed from the operation
tables. The construction need not produce a small basis.

\subsection{Catalan varieties}

For $r\ge2$, let $[r]=\{1,\ldots,r\}$ with its usual order, and let
$\Cminus r$ consist of all order-preserving decreasing maps
$[r]\to[r]$. Addition is pointwise maximum, and multiplication is
composition, with maps acting on the right:
$i(\alpha\beta)=(i\alpha)\beta$.
The resulting semiring and its identities were studied by Gusev and
Volkov~\cite{GusevVolkovCatalan}.

The maximum of two such maps is again order-preserving and decreasing.
Composition distributes over maximum on both sides: on one side this
is pointwise evaluation, and on the other it follows because an
order-preserving map of a chain preserves binary maximum. Thus
$\Cminus r$ is an ai-semiring. Further,
\[
 i(\alpha\beta)\le i\alpha,
 \qquad
 i(\alpha\beta)=(i\alpha)\beta\le i\beta.
\]
The first inequality uses that $\beta$ is decreasing, and the second
uses that $\alpha$ is decreasing and $\beta$ is order-preserving.
Therefore $\Cminus r$ is an incline.

As recalled in the introduction, Gusev and
Volkov~\cite[Theorem~4.12]{GusevVolkovCatalan} proved
\begin{equation}\label{eq:catalan}
 \var\Cminus r=\W_r\qquad(r\ge3).
\end{equation}
Theorem~\ref{thm:uniform} gives the following strengthening of their
finite-basis theorem.

\begin{corollary}\label{cor:catalan}
For every $r\ge2$, the variety $\var\Cminus r$ is a Cross variety.
For $r\ge3$, every subvariety is determined, within $\W_r$, by its
identities on at most $r-2+2^{r-1}$ variables. In particular, the
subvariety lattice of $\var\Cminus3$ is finite, and each of its
members is determined within $\W_3$ by identities on at most five
variables.
\end{corollary}

\begin{proof}
The first assertion follows from Theorem~\ref{thm:main}, since
$\Cminus r$ is finite. For $r\ge3$, use~\eqref{eq:catalan} and
Theorem~\ref{thm:uniform}. Finally, $K_3=1+2^2=5$.
\end{proof}

The inclusion $\Cminus r\in\W_r$ for $r\ge3$ can also be verified
directly. Evaluate a product of $r-1$
decreasing maps at a point of $[r]$. If any step is stationary,
repeating that factor $r$ times leaves the final value unchanged.
If every step is strict, the initial point must be $r$ and the final
point must be $1$. Every expanded product is below the original
product by contraction, so all expanded products also take the value
$1$ at that point. In either case, the pointwise maximum of the
expanded products equals the original product. This proves $\sigma_r$.

\begin{remark}
Theorem~\ref{thm:main} concerns the varieties generated by individual
finite inclines. The whole variety of inclines is finitely axiomatized
by $\AI$ and~\eqref{eq:contraction}, and is generated by all the
decreasing Catalan semirings together
\cite[Theorem~5.11]{GusevVolkovCatalan}. It is not locally finite.
For example, the positive integers with
$a+b:=\min\{a,b\}$ and $a\cdot b:=a+b$ (ordinary addition on the
right) form an incline generated by the integer $1$, and this incline
is infinite. Thus the local-finiteness assertion of
Theorem~\ref{thm:uniform} depends on fixing a stabilization identity
$\sigma_n$. The separate question in
\cite[Section~6.3]{GusevVolkovCatalan} about the cardinality of the
subvariety lattice of the whole variety of inclines is not settled by
these results.
\end{remark}


\end{document}